\documentclass[11pt]{amsart}

\usepackage{amsmath,amssymb,amscd}
\usepackage{xcolor}

\usepackage{amsmath,amssymb}

\usepackage{amsxtra, amsmath}
\usepackage{amssymb, amscd}
\usepackage{graphicx}
\usepackage{mathrsfs}

\usepackage{mathtools}

\newtheorem{thm}{Theorem}[section]
\newtheorem{prop}[thm]{Proposition}

\theoremstyle{definition}
\newtheorem{definition}[thm]{Definition}

\theoremstyle{remark}
\newtheorem{remark}[thm]{Remark}

\usepackage{dsfont}
\allowdisplaybreaks

\title{Singular examples of the matrix  Bochner problem}
\author{Ignacio Bono Parisi}
\author{Ines Pacharoni}



\address{CIEM-FaMAF\\ Universidad Nacional de C\'or\-do\-ba\\
CP 5000, C\'or\-do\-ba,  Argentina}
\email{ignacio.bono@unc.edu.ar, ines.pacharoni@unc.edu.ar}

\begin{document}
\begin{abstract}
The Matrix Bochner Problem aims to classify which weight matrices have their sequence of orthogonal polynomials as eigenfunctions of a second-order differential operator. Casper and Yakimov, in \cite{CY18}, demonstrated that, under certain hypotheses, all solutions to the Matrix Bochner Problem are noncommutative bispectral Darboux transformations of a direct sum of classical scalar weights. This paper aims to provide the first proof that there are solutions to the Matrix Bochner Problem that do not arise through a noncommutative bispectral Darboux transformation of any direct sum of classical scalar weights. This initial example could contribute to a more comprehensive understanding of the general solution to the Matrix Bochner Problem. 
\end{abstract} 
\maketitle

\section{Introduction}

Back in 1929, Bochner \cite{B29} posed the problem of determining all scalar-valued orthogonal polynomials that are eigenfunctions of some arbitrary, but fixed, second-order differential operator.
Bochner proved that, up to an affine change of coordinates, the only weights satisfying these properties are the classical weights $e^{-x^{2}}$, $e^{-x}x^{\alpha}$,
and $(1-x)^\alpha(1 + x)^\beta$ 
of Hermite, Laguerre, and Jacobi respectively. 

Orthogonal matrix polynomials are sequences of matrix-valued polynomials that are pairwise orthogonal with respect to a matrix-valued inner product defined by an $N\times N$ weight matrix $W(x)$.
The theory of these matrix-valued orthogonal polynomials, without any consideration of differential equations, goes back to \cite{K49} and \cite{K71}. In \cite{D97}, the study of matrix-valued orthogonal polynomials that are eigenfunctions of certain second-order symmetric differential operators was initiated.

\smallskip

Nowadays, the problem of finding weight matrices $W(x)$ of size $N \times N$  such that the associated sequence of orthogonal matrix polynomials are eigenfunctions of a second-order matrix differential operator is known as the Matrix Bochner Problem.
  
\smallskip

The first nontrivial solutions of the Matrix Bochner Problem were given in  \cite{DG04}, \cite{G03}, \cite{GPT01}, and \cite{GPT02}, by using different methods to find them. Since then, over the past twenty years, several additional examples have been discovered. 
See \cite{GPT05}, \cite{CG06}, \cite{DG07}, \cite{PT07}, \cite{DdI08}, \cite{DdI08b}, \cite{DG05} \cite{CMV05}, \cite{CMV07}, \cite{P08}, \cite{PR08}, \cite{PZ16}, \cite{KPR12}, \cite{KRR14}.  

\smallskip

 In \cite{CY18}, Casper and Yakimov developed a general framework for the study of the Matrix Bochner Problem, performing an important breakthrough in this area. They proved that, under certain hypotheses,  every weight matrix $W(x)$ can be obtained from classical scalar weights by a noncommutative bispectral Darboux transformation of a diagonal weight matrix whose entries are classical Jacobi, Laguerre, or Hermite weights. 
One of the main theorems in \cite{CY18} is the following classification result.
\begin{thm}\label{clas1}
(\cite{CY18}, Theorem 1.3). Let $W(x)$ be a weight matrix and suppose that $\mathcal{D}(W)$ contains a $W$-symmetric second-order differential operator
$D = \partial^{2}G_{2}(x) + \partial G_{1}(x) + G_{0}(x)$,
with $G_{2}(x)W(x)$ positive-definite on the support of $W(x)$. 
Then  $W(x)$ is a noncommutative bispectral Darboux transformation of a direct sum of classical weights if and only if the algebra $\mathcal{D}(W)$ is full.
\end{thm}

The algebra $\mathcal D(W)$ is the set of all differential operators having any sequence of matrix-valued orthogonal polynomials with respect to $W$ as eigenfunctions. 
The algebraic condition for this algebra to be considered full is somewhat technical and it means that there exist nonzero $W$-symmetric operators $\mathfrak{D}_{1},\ldots,\mathfrak{D}_{N}$ in $\mathcal{D}(W)$,  such that  $\mathfrak{D}_{i}\mathfrak{D}_{j} = 0 $,  for $i\not= j\quad $ with $ \mathfrak{D}_{1} + \cdots + \mathfrak{D}_{N}$ a central element in $\mathcal{D}(W)$ which is not a zero divisor. 

This result shows that the algebraic structure of the algebra $\mathcal D(W)$ has a profound influence on the shape of the weight matrix $W(x)$ itself.
 In the classical scalar cases of Hermite, Laguerre, and Jacobi weights, the structure of this algebra is well understood: it is a polynomial algebra in a second-order differential operator.
In the matrix case, the basic definitions and main results concerning the algebra $\mathcal D(W)$ were given in \cite{GT07}. 
However, understanding the structure of this algebra, even in relatively simple $2\times 2$ examples, is certainly a non-trivial problem. See for example \cite{T11},  \cite{Z15}.

\smallskip

Theorem \ref{clas1} can be seen as a solution to the Matrix Bochner Problem under 
the natural conditions that $\mathcal D(W)$ is full and the leading coefficient of the second-order differential operator multiplied by the weight $W(x)$ is positive definite. However, there exist examples that escape from this general framework. This is the case, for instance, of the following Hermite-type weights
\begin{equation}\label{weight}
     W(x) = e^{-x^{2}}\begin{pmatrix} e^{2bx} + a^{2}x^{2} && ax \\
ax && 1 \end{pmatrix}, \quad   \text { $a,b,  \in \mathbb{R}$,  $a,b \neq 0$ }, \quad x \in \mathbb{R}, 
\end{equation}
introduced first in \cite{CG06}, for $a=b=1$. 
The differential operator 
\begin{equation} \label{operatorD}
D = \partial^{2}I + \partial \begin{pmatrix} -2x + 2b && -2abx + 2a \\
0 && -2x\end{pmatrix} + \begin{pmatrix} -2 && 0 \\ 0 && 0 \end{pmatrix}. 
\end{equation}
is a $W$-symmetric operator in the algebra $\mathcal D(W)$ and hence the weight $W$ is a solution to the Matrix Bochner Problem. In this paper, we prove the following result:

\begin{thm}\label{not Darboux}
The weight matrix $W(x)$ is not a noncommutative bispectral Darboux transformation of any direct sum of classical scalar weights.
\end{thm}

\smallskip

The paper is organized as follows. Section \ref{Backgr} contains the background material on matrix-valued orthogonal polynomials with respect to a weight function $W$, the $W$-symmetric operators in the algebra $\mathcal D(W)$ and Fourier algebras associated to a weight $W$.
In Section \ref{structure}, we ascertain the algebraic structure of $\mathcal{D}(W)$. Initially, we investigate the centralizer of $D$ in $\mathcal{D}(W)$ and demonstrate that it is a polynomial algebra over $\mathbb{C}$ in the operator $D$.
Then, by using that $\mathcal D(W)$ has a nontrivial center, we prove that any differential operator in $\mathcal D(W)$ commutes with $D$, and thus, it is a polynomial in $D$.  In this case, it is easy to verify that the algebra $\mathcal D (W)$ is not full, and therefore, it is not a noncommutative bispectral Darboux transformation of any direct sum of classical scalar weights. 

In Section \ref{RFA}, we explicitly describe the right Fourier algebra $\mathcal{F}_{R}(W)$ of $W$ by leveraging the following relationship between $W$ and the direct sum of scalar Hermite weights $\widetilde{W}=\left(\begin{smallmatrix} e^{-x^2+2bx} && 0 \\ 0 && e^{-x^2} \end{smallmatrix}\right)$,
\begin{equation*}
W(x) = T(x) \widetilde{W}(x) T(x)^* ,\qquad \text { for some matrix polynomial $T(x)$. }
\end{equation*}

 We provide the explicit expressions of a sequence of orthogonal polynomials for $W$ and the three-term recursion relation satisfied by them in  Section \ref{mop}. Finally, in Section \ref{HDEx}, we present another case of a weight matrix of size $3\times 3$ with similar properties. It seems reasonable to assume that such examples can be found for matrices of any size.
A natural question arising from this first example concerns its level of exceptionality. How many examples of this kind exist, and what insights can we derive from them?
Moreover, can similar examples be found in different dimensions or contexts related to Laguerre or Jacobi scalar weights?
We hope that this first example could contribute to a more comprehensive understanding of the general solution to the Matrix Bochner Problem.

\section{Background}\label{Backgr}

\subsection{Matrix-valued orthogonal polynomials and  the algebra $\mathcal D(W)$}\label{sec-MOP}

Let $W=W(x)$ be a weight matrix of size $N$ on the real line, that is, a complex $N\times N$ matrix-valued integrable function on the interval $(x_0,x_1)$ such that $W(x)$ is positive definite almost everywhere and with finite moments of all orders. Let $\operatorname{Mat}_N(\mathbb{C})$ be the algebra of all $N\times N$ complex matrices and let $\operatorname{Mat}_N(\mathbb{C}[x])$ be the algebra of polynomials in the indeterminate $x$ with coefficients in $\operatorname{Mat}_N(\mathbb{C})$. We consider the following Hermitian sesquilinear form in the linear space $\operatorname{Mat}_N(\mathbb{C}[x])$
\begin{equation*}
  \langle P,Q \rangle =  \langle P,Q \rangle_W = \int_{x_0}^{x_1} P(x) W(x) Q(x)^*\,dx.
\end{equation*}

Given a weight matrix $W$ one can construct sequences $\{Q_n\}_{n\in\mathbb{N}_0}$ of matrix-valued orthogonal polynomials, i.e. the $Q_n$ are polynomials of degree $n$ with nonsingular leading coefficient and $\langle Q_n,Q_m\rangle=0$ for $n\neq m$.
We observe that there exists a unique sequence of monic orthogonal polynomials $\{P_n\}_{n\in\mathbb{N}_0}$ in $\operatorname{Mat}_N(\mathbb{C}[x])$.
By following a standard argument (see \cite{K49} or \cite{K71}) one shows that the monic orthogonal polynomials $\{P_n\}_{n\in\mathbb{N}_0}$ satisfy a three-term recursion relation
\begin{equation}\label{ttrr}
    x P_n(x)=P_{n+1}(x) + B_{n}P_{n}(x)+ C_nP_{n-1}(x), \qquad n\in\mathbb{N}_0,
\end{equation}
where $P_{-1}=0$ and $B_n, C_n$ are matrices depending on $n$ and not on $x$.

Along this paper, we consider that an arbitrary matrix differential operator
\begin{equation}\label{D2}
  {D}=\sum_{i=0}^s \partial ^i F_i(x),\qquad \partial=\frac{d}{dx},
\end{equation}
acts on the right on a matrix-valued function $P$ i.e.  
\begin{equation}\label{actionD}
    (P{D})(x)=\sum_{i=0}^s \partial ^i (P)(x)F_i(x).
\end{equation}
\noindent 
We consider the algebra of these operators with polynomial coefficients 
$$\operatorname{Mat}_{N}(\Omega[x])=\Big\{D = \sum_{j=0}^{n} \partial^{j}F_{j}(x) \, : F_{j} \in \operatorname{Mat}_{N}(\mathbb{C}[x]) \Big \}.$$
\noindent 
More generally, when necessary, we will also consider $\operatorname{Mat}_{N}(\Omega[[x]])$, the set of all differential operators with coefficients in $\mathbb{C}[[x]]$, the ring of power series with coefficients in $\mathbb{C}$.

\begin{prop}[\cite{GT07}, Propositions 2.6 and 2.7]\label{eigenvalue-prop}
 Let $W=W(x)$ be a weight matrix of size $N\times  N$, and let $\{P_n\}_{n\geq 0}$ be the sequence of monic orthogonal polynomials in $\operatorname{Mat}_N(\mathbb{C}[x])$. If $D$ is a differential operator of order $s$, as in \eqref{D2}, such that
$$P_nD=\Lambda_n P_n, \qquad \text{for all } n\in\mathbb{N}_0,$$
with $\Lambda_n\in \operatorname{Mat}_N(\mathbb{C})$, then
$F_i(x)=\sum_{j=0}^i x^j F_j^i$, $F_j^i \in \operatorname{Mat}_N(\mathbb{C})$, is a polynomial, and $\deg(F_i)\leq i$. Moreover, $D$ is determined by the sequence $\{\Lambda_n\}_{n\geq 0}$, and
\begin{equation}\label{eigenvaluemonicos}
   \Lambda_n=\sum_{i=0}^s [n]_i F_i^i, \qquad \text{for all } n\geq 0,
\end{equation}
where $[n]_i=n(n-1)\cdots (n-i+1)$, $[n]_0=1$.
\end{prop}

Given a weight matrix $W$, the algebra 
\begin{equation}\label{algDW}
  \mathcal D(W)=\{D\in \operatorname{Mat}_{N}(\Omega[x])\, : \, P_nD=\Lambda_n(D) P_n, \, \Lambda_n(D)\in \operatorname{Mat}_N(\mathbb{C}), \text{ for all }n\in\mathbb{N}_0\}
\end{equation}
is introduced in \cite{GT07}, where $\{P_n\}_{n\in \mathbb{N}_0}$ is any sequence of matrix-valued orthogonal polynomials with respect to $W$.

We observe that the definition of $\mathcal D(W)$ depends only on the weight matrix $W$ and not on the particular sequence of orthogonal polynomials. This is because two sequences $\{P_n\}_{n\in\mathbb{N}_0}$ and $\{Q_n\}_{n\in\mathbb{N}_0}$ of matrix orthogonal polynomials with respect to the weight $W$ are related by $P_n=M_nQ_n$, where $\{M_n\}_{n\in\mathbb{N}_0}$ consists of invertible matrices (see \cite[Corollary 2.5]{GT07}).

\begin{prop} [\cite{GT07}, Proposition 2.8]\label{prop2.8-GT}
For each $n\in\mathbb{N}_0$, the mapping $D\mapsto \Lambda_n(D) $ is a representation of $\mathcal D(W)$ in $\operatorname{Mat}_N(\mathbb{C})$. Moreover, the sequence of representations $\{\Lambda_n\}_{n\in\mathbb{N}_0}$ separates the elements of $\mathcal D(W)$.
\end{prop}

The {\em formal adjoint} on $\operatorname{Mat}_{N}(\Omega([[x]])$, denoted  by $\mbox{}^*$, is the unique involution extending Hermitian conjugate on $\operatorname{Mat}_{N}(\mathbb C)[x]$ and sending $\partial I$ to $-\partial I$. 
The {\em formal $W$-adjoint} of $\mathfrak{D}\in \operatorname{Mat}_{N}(\Omega[x])$, is the differential operator $\mathfrak{D}^{\dagger} \in \operatorname{Mat}_{N}(\Omega[[x]])$ defined by
$$\mathfrak{D}^{\dagger}:= W(x)\mathfrak{D}^{\ast}W(x)^{-1},$$
where $\mathfrak{D}^{\ast}$ is the formal adjoint of 
$\mathfrak D$. 
An operator $\mathfrak{D}\in \operatorname{Mat}_{N}(\Omega[x])$ is called {\em $W$-adjointable} if there exists $\widetilde {\mathfrak{D}} \in \operatorname{Mat}_{N}(\Omega[x])$, such that
$$\langle P\mathfrak{D},Q\rangle=\langle P,Q\widetilde{\mathfrak{D}}\rangle,$$ for all $P,Q\in \operatorname{Mat}_N(\mathbb{C}[x])$. Then we say that the operator $\widetilde{\mathfrak D}$ is the $W$-adjoint of $\mathfrak D $.

\begin{prop}[\cite{CY18}, Prop. 2.23] \label{adjuntas}  
If 
 $\mathfrak{D} \in \operatorname{Mat}_{N}(\Omega[x])$ is $W$-adjointable and $\mathfrak{D}^{\dagger} \in \operatorname{Mat}_{N}(\Omega[x])$, then $\mathfrak{D}^{\dagger}$ is the $W$-adjoint of $\mathfrak D$, i.e.
 $$\langle \, P\mathfrak{D},Q\, \rangle=\langle \, P,\, Q {\mathfrak{D}}^\dagger \rangle,$$ for all $P,Q\in \operatorname{Mat}_N(\mathbb{C}[x])$.
 \end{prop}
 
For $\mathfrak D= \sum_{j=0}^n \partial ^j F_j \in \operatorname{Mat}_{N}(\Omega[x]) $,  the formal $W$-adjoint of $\mathfrak D$ is given by $\mathfrak{D}^\dagger= \sum_{k=0}^n \partial ^k G_k$, with
\begin{equation}\label{daga}
    G_k=  \sum_{j = 0}^{n-k}(-1)^{n-j} \binom{n-j}{k}
    (WF_{n-j}^*)^{(n-k-j)} W^{-1} , \qquad \text{for } 0\leq k\leq n.  
 \end{equation}
 It is a matter of careful integration by parts to see that $\langle \, P\mathfrak{D}, Q\, \rangle=\langle \, P,\, Q {\mathfrak{D}}^\dagger \rangle$ if  the following set of ``boundary conditions" are satisfied
  \begin{equation} 
  \lim_{x\to x_i}\, \sum_{j = 0}^{p-1}(-1)^{n-j+p-1} \binom{n-j}{k}
    \big (F_{n-j}(x)W(x)\big)^{(p-1-j)}=0, 
      \end{equation}
      for $1\leq p\leq n$  and $0\leq k\leq n-p$,
      where $x_i$ are the endpoints     
      of the support of the weight $W$.

\

We say that a differential operator $D\in \mathcal D(W)$ is $W$-{\em symmetric} if $\langle P{D},Q\rangle=\langle P,QD\rangle$, for all $P,Q\in \operatorname{Mat}_N(\mathbb{C}[x])$.
An operator $\mathfrak{D}\in \operatorname{Mat}_{N}(\Omega[x])$ is called {\em formally } $W${\em -symmetric} if $\mathfrak{D}^{\dagger} = \mathfrak{D}$.
In particular if $\mathfrak{D}\in \mathcal{D}(W)$, then $\mathfrak{D}$ is $W$-symmetric if and only if it is formally $W$-symmetric.

It is shown that the set $\mathcal S(W)$ of all symmetric operators in $\mathcal D(W)$ is a real form of the space $\mathcal D(W)$, i.e.
$$\mathcal D(W)= \mathcal S (W)\oplus i \mathcal S (W),$$ as real vector spaces. 

The condition of symmetry for a differential operator in the algebra $\mathcal D(W)$ is equivalent to the following set of differential equations involving the weight $W$ and the coefficients of $D$. (cf. \cite{dI07}, Theorem 2.1.3).

\begin{thm}\label{equivDsymm}
 Let $\mathfrak{D} =\sum_{i=0}^n \partial^i F_i(x)$ be a  differential operator of order $n$ in $\mathcal{D}(W)$.
 Then $\mathfrak{D}$ is $W$-symmetric if and only if
\begin{equation*}
    \sum_{j = 0}^{n-k}(-1)^{n-j} \binom{n-j}{k}
    (F_{n-j}W)^{(n-k-j)}= WF_{k}^\ast   
\end{equation*}
for all $0\leq k \leq n$.
 \end{thm}
 \begin{proof}  
 An operator $\mathfrak{D}\in \mathcal{D}(W)$ is $W$-symmetric if and only if $\mathfrak{D}= \mathfrak{D}^\dagger$. The statement follows by using the explicit expression of the coefficients of  $\mathfrak{D}^\dagger$ given in  \eqref{daga}. 
 \end{proof}

In particular, the coefficients of a differential operator of order two in $\mathcal S(W)$ satisfy the classical symmetry equations obtained in \cite{GPT03a}.
\begin{align*}
  F_2 W &=WF_2^*\\
    2(F_2W)'&-F_1W =WF_1^*\\
    (F_2W)''&-(F_1W)'+F_0W=WF_0^*
\end{align*}

\smallskip

\begin{definition}\label{D(W)full}
    Let $W$ a $N\times N$ weight matrix function. We say that the algebra $\mathcal{D}(W)$ is full if there exist nonzero $W$-symmetric operators $\mathfrak{D}_{1},\ldots,\mathfrak{D}_{N}$ in $\mathcal{D}(W)$, such that  $\mathfrak{D}_{i}\mathfrak{D}_{j} = 0 $,  for $i\not= j\quad $ with $ \mathfrak{D}_{1} + \cdots + \mathfrak{D}_{N}$ a central element in $\mathcal{D}(W)$ which is not a zero divisor. 
\end{definition}

\subsection{Bispectral Darboux transformations}
We recall the notion of right Fourier algebra associated with a weight matrix $W(x)$ as given in \cite{CY18}. 
We consider the space of all semi-infinite sequences of matrix-valued rational functions 
$$ \mathcal{P} =\{ P:\mathbb{C}\times \mathbb{N}_0 \longrightarrow M_N(\mathbb{C}) \, : \, P(x,n)  \text{ is a rational function of $x$, for each fixed $n$}  \}.$$ 
We denote by $\operatorname{Mat}_{N}(\mathcal{S})$ the algebra of all discrete operators of the form 
 $\mathscr{M}= \sum_{j=-\ell}^{k} A_{j}(n)\, \delta^j$, $A_{j}(n)\in\operatorname{Mat}_{N}(\mathbb{C})$, with its left action defined as
\begin{equation}\label{discreteop}(\mathscr{M}\cdot P)(x,n)= \sum_{j=-\ell}^{k} A_{j}(n)(\delta^j \cdot P)(x,n)= \sum_{j=-\ell}^k A_j(n) P(x,n+j).
\end{equation}
 
The right Fourier algebra associated with a weight matrix $W$ (or associated with its (unique) sequence of monic orthogonal polynomials) is defined by 
\begin{equation*}
\begin{split}
    \mathcal F_R(W)= \mathcal F_R(P) &=\{ \mathfrak D\in \operatorname{Mat}_{N}(\Omega[x])\, : \, \exists \, \mathscr{M} \in \operatorname{Mat}_{N}(\mathcal{S})\text{ such that } P(x,n)\cdot \mathfrak D=\mathscr{M}\cdot P(x,n) \}.
\end{split}
\end{equation*}

\smallskip

In \cite{CY18}, Theorem 3.7, the authors give an explicit description of the right Fourier algebra associated with a weight matrix $W$. 
\begin{equation} \label{fourier algebra}
\mathcal{F}_{R}(W) = \{ \mathfrak{D}\in \operatorname{Mat}_{N}(\Omega[x]):\mathfrak{D} \text{ is } W\text{-adjointable and } \mathfrak{D}^{\dagger}\in \operatorname{Mat}_{N}(\Omega[x]) \}.
\end{equation}

\begin{remark}\label{adj}
    In the case of the weight matrices considered in this paper, if
$\mathfrak{D}\in \operatorname{Mat}_{N}(\Omega[x])$ such that   $\mathfrak{D}^{\dagger}\in \operatorname{Mat}_{N} (\Omega[x])$, then $\mathfrak D$ is $W$-adjointable due to the exponential decay of the weight at $\pm \infty$.
\end{remark}

\begin{definition}\label{Darbou-transf-def}
Let $W(x)$ and $\widetilde{W}(x)$ be weight matrices, and let $P(x,n)$ and $\widetilde{P}(x,n)$ be their associated sequences of monic orthogonal polynomials. We say that $\widetilde{P}(x,n)$ is a bispectral Darboux transformation of $P(x,n)$ if there exist differential operators $\mathfrak{D},\widetilde{\mathfrak{D}} \in \mathcal{F}_{R}(W)$, polynomials $F(x),\widetilde{F}(x)$, and sequences of matrices $C(n),\widetilde{C}(n)$ which are nonsingular for almost every $n$ and satisfy
$$C(n)\widetilde{P}(x,n) = P(x,n) \cdot \mathfrak{D}F(x)^{-1} \text{ and } \quad  \widetilde{C}(n)P(x,n) = \widetilde{P}(x,n) \cdot \widetilde{F}(x)^{-1}\widetilde{\mathfrak{D}}.$$
We say that $\widetilde{W}(x)$ is a {\em noncommutative bispectral Darboux transformation } of $W(x)$ if $\widetilde{P}(x,n)$ is a bispectral Darboux transformation of $P(x,n)$.
\end{definition}

\section{The structure of the algebra $\mathcal{D}(W)$}\label{structure}

This section aims to prove that the algebra $\mathcal D(W)$, for  the Hermite weight  
$$W(x) = e^{-x^{2}}\begin{pmatrix} e^{2bx} + a^{2}x^{2} && ax \\
ax && 1 \end{pmatrix}, \qquad  
a,b \neq  0 ,  \quad 
x \in \mathbb{R}
$$
is a polynomial algebra on the $W$-symmetric differential operator
$$D = \partial^{2}I + \partial \begin{pmatrix} -2x + 2b && -2abx + 2a \\
0 && -2x\end{pmatrix} + \begin{pmatrix} -2 && 0 \\ 0 && 0 \end{pmatrix}.$$

 We will first prove that the centralizer of $D$ in $\mathcal D(W)$,
$$\mathcal Z_{\mathcal D(W)}(D)=\{\mathfrak{D}\in \mathcal D(W): \mathfrak D D=D\mathfrak D\},$$
is a polynomial algebra in $D$. Subsequently, we will demonstrate that any differential operator in $\mathcal D(W)$ commutes with $D$.

\ 
 
A differential operator $\mathfrak D\in \mathcal D(W)$ is of the form $\mathfrak D= \displaystyle \sum_{j=0}^{n}\frac{d^j}{dx^{j}} \,F_j$, where $F_j=F_j(x)$ are polynomial matrices of degree at most $j$. First of all, we shall prove that all these coefficients $F_j$ are upper triangular matrices.
 
 \begin{prop}\label{(2,1)} 
The coefficients of any operator $\mathfrak{D} \in \mathcal{D}(W)$ are upper triangular matrices. 
 \end{prop}
 \begin{proof} We can assume that $\mathfrak D$ is a symmetric operator because   $\mathcal{D}(W) = \mathcal{S}(W) \oplus i\mathcal{S}(W)$. 
 
  Let $\mathfrak{D} = \sum_{j=0}^{n}\partial^{j}F_{j} \in \mathcal{S}(W)$ with $F_{j} = \begin{pmatrix} p_{j} && r_{j} \\ g_{j} && q_{j} \end{pmatrix}$, and  $p_{j},r_{j},g_{j},q_{j} \in \mathbb{C}[x]$, 
 for all $0 \leq j \leq n$. 
   From Theorem \ref{equivDsymm}, the coefficients of $\mathfrak{D}$ satisfy the following set of differential equations, for $k=0,\dots, n$,
   \begin{equation} \label{aux}
   \sum_{j = 0}^{n-k}(-1)^{n-j} \binom{n-j }{k}(F_{n-j}W)^{(n-k-j)}= WF_{k}^\ast.
   \end{equation}
For each $0\leq k \leq n$, the entry $(1,2)$ in \eqref{aux} gives  
\begin{equation*}
    e^{-x^{2}}(a^{2}x^{2}\overline{g}_{k} + ax \overline{q}_{k}) + e^{-x^{2}+2bx} \overline{g}_{k} = \sum_{j = 0}^{n-k}(-1)^{n-j} \binom{n-j }{k}\big (e^{-x^{2}}(ax p_{n-j} + r_{n-j})\big)^{(n-k-j)}.
\end{equation*}
Hence we obtain that $e^{2bx}\bar g_k$ is a polynomial function and therefore $g_k$  must be zero. 
 \end{proof}

Now, we characterize the equations that the coefficients of a symmetric operator in the algebra $\mathcal{D}(W)$ must satisfy.

\begin{prop}
A differential operator  $\mathfrak{D} = \sum_{j = 0}^{n} \partial^{j}F_{j} \in
\mathcal{D}(W)$ is  symmetric if and only if its polynomial coefficients $F_{j} = \begin{pmatrix} p_{j} && r_{j} \\ 0 && q_{j} \end{pmatrix}$ satisfy the following  set of equations, for each $0\leq k\leq n$

\begin{equation}\label{A4}
     \sum_{j=0}^{n-k-1}(-1)^{n-j} \binom{n-j }{k} \left(e^{-x^{2}}q_{n-j}\right)^{(n-k-j)} = e^{-x^{2}}(\overline{q}_{k}+(-1)^{k+1}q_{k}),
\end{equation}
\begin{equation}\label{A2}
    \sum_{j=0}^{n-k-1}(-1)^{n-j} \binom{n-j }{k} \left(e^{-x^{2}}(axp_{n-j} + r_{n-j})\right)^{(n-k-j)} = e^{-x^{2}}(ax\overline{q}_{k}+(-1)^{k+1}(axp_{k}+r_{k})),
\end{equation}
\begin{equation}\label{A3}
     \sum_{j=0}^{n-k-1}(-1)^{n-j} \binom{n-j }{k} a(n-k-j)\left (e^{-x^{2}}q_{n-j}\right )^{(n-k-j-1)} = e^{-x^{2}}(
     ax(\overline{p}_{k}-\overline{q}_{k}) + \overline{r}_{k} ).
\end{equation}
\begin{equation}\label{A1}
    \begin{split}
        \sum_{j=0}^{n-k-1}(-1)^{n-j} \binom{n-j }{k}& \left(e^{-x^{2}} (a^{2}x^{2}p_{n-j}+axr_{n-j}\right ))^{(n-k-j)}  \\ 
        & = e^{-x^{2}}a^{2}x^{2}(\overline{p}_{k}+(-1)^{k+1}p_{k}) 
         \quad + e^{-x^{2}}ax(\overline{r}_{k}+(-1)^{k+1}r_{k}),
    \end{split}
\end{equation}
\begin{equation}\label{AA1}
     \sum_{j=0}^{n-k-1}(-1)^{n-j} \binom{n-j }{k} \left(e^{-x^{2}+2bx}p_{n-j}\right)^{(n-k-j)} = e^{-x^{2}+2bx}(\overline{p}_{k}+(-1)^{k+1}p_{k}).
\end{equation}
\end{prop}
\begin{remark}
For \( k = n \), the left-hand side of equations \eqref{A4}-\eqref{AA1} must be taken as zero. Additionally, we adopt the convention that if \( n < m \), the sum \(\sum_{j=m}^n a_j\) is equal to zero.
\end{remark}

\begin{proof}
From Theorem \ref{equivDsymm} we have that the coefficients of $\mathfrak{D}$  satisfy 
\begin{equation*}
    \sum_{j=0}^{n-k-1}(-1)^{n-j} \binom{n-j }{k}(F_{n-j}W)^{(n-k-j)} = WF_{k}^{\ast} + (-1)^{k+1}F_{k}W, \qquad \text { for }k=0,\dots, n.
\end{equation*}

The entries (2,2) and (1,2) in  
the above matrix equation give the equations \eqref{A4} and \eqref{A2}. The entry (2,1) is 
\begin{equation}\label{Aa3}
    \sum_{j=0}^{n-k-1}(-1)^{n-j} \binom{n-j }{k} \left(e^{-x^{2}}axq_{n-j}\right)^{(n-k-j)} = e^{-x^{2}}(ax\overline{p}_{k}+\overline{r}_{k}+(-1)^{k+1}axq_{k}).
\end{equation}
By using that
\begin{equation*}
        (e^{-x^{2}}axq_{n-j})^{(n-k-j)} 
        = ax (e^{-x^{2}}q_{n-j})^{(n-k-j)}+a(n-k-j)(e^{-x^{2}}q_{n-j})^{(n-k-j-1)},
\end{equation*}
and combining with 
\eqref{A4} 
we obtain that \eqref{Aa3} is equivalent to \eqref{A3}.
  Finally, the entry (1,1) is the sum of equations \eqref{A1} and \eqref{AA1}, which splits into two equations due to the factor $e^{2bx}$. This concludes the proof.
\end{proof}

Observe that for a symmetric operator $\mathfrak D $, the identity  \eqref{A3} gives  an expression of the polynomial $r_k$ in terms of the coefficients  $p_j$ and $q_j$, for 
$0 \leq k \leq n$
\begin{equation}\label{r}
    r_{k} = ax(q_{k}-p_{k}) + \sum_{j=0}^{n-k-1}(-1)^{n-j} \binom{n-j }{k} a(n-k-j)(e^{-x^{2}}\overline{q}_{n-j})^{(n-k-j-1)}e^{x^{2}}.
\end{equation}

\

The differential operator $D$, can be written as  $D = \partial^{2}I + \partial (Ax+B) + C$, where
$$A = \begin{pmatrix} -2 && -2ab \\ 0 && -2 \end{pmatrix}, \quad B = \begin{pmatrix} 2b && 2a \\ 0 && 0 \end{pmatrix}, \text{ and } \quad C = \begin{pmatrix} -2 && 0 \\ 0 && 0 \end{pmatrix}.$$

Let $\mathfrak{D}= \sum_{j=0}^{n}\partial^j  F_{j} \in \mathcal{D}(W)$ be a differential operator of order $n$. By Proposition \ref{(2,1)}, we can write the coefficients  
$$F_{j} = \begin{pmatrix} p_{j} & r_{j} \\ 0 & q_{j} \end{pmatrix},$$  $p,q,r \in \mathbb{C}[x]$  polynomials with degree $\leq j$.
Therefore, 
$\mathfrak{D}$ commutes with  $D$  if and only if the matrix coefficients satisfy
\begin{equation}\label{dD=Dd}
    2F_{k-1}'+F_{k-1}(Ax+B)-(Ax+B)F_{k-1}
    =kAF_{k}+CF_{k}-F_{k}C-F_{k}''-F_{k}'(Ax+B),
\end{equation}

\noindent for all $k$. 
As usual, we assume that 
$F_{j} = 0$ for all $j \notin \{0,1,\cdots,n\}$.

\smallskip
The equation \eqref{dD=Dd}, for $k=n+1$, gives  
\begin{equation}\label{dD=Ddn}
    2F_{n}' + F_{n}(Ax+B) -(Ax+B)F_{n}=0.
\end{equation}
Thus, we get 
\begin{equation*}
    \begin{pmatrix}
    2p_n'& 2r_n'\\ 0& 2q_n'
    \end{pmatrix} = 
    \begin{pmatrix}
    0&(2abx-2a)(q_n-p_n)+2br_n\\0&0
    \end{pmatrix}.
\end{equation*}
From here, we obtain that there exist $\alpha,\beta \in \mathbb{C}$ such that  
\begin{equation} \label{k=0} 
p_n=\alpha, \qquad  q_n=\beta , \quad \text{ and } 
\quad r_n=a(\alpha-\beta) x
\end{equation}

\noindent Again  from \eqref{dD=Dd}, now with $k=n$, we have that 
\begin{equation*}
    2F_{n-1}' + F_{n-1}(Ax+B)- (Ax+B)F_{n-1} =  CF_{n}-F_{n}C + nAF_{n}-F_{n}'(Ax+B)-F_{n}''.
\end{equation*}
By looking at the entries $(1,1)$ and $(2,2)$ in the above equation,  we see that
$p_{n-1}' = -n\alpha $ and  $q_{n-1}' = -n\beta$.
Then
\begin{equation}\label{k=1}
        p_{n-1}  = -n\alpha x + \alpha_{1}, \qquad 
        q_{n-1}  = -n\beta x + \beta_{1}
  \end{equation}
with $\alpha_{1},\beta_{1} \in \mathbb{C}$.

\medskip
We introduce the following notation: 
$[\![n]\!]_{k} = n(n-2)\cdots(n-2(k-1))$, $[\![n]\!]_{0} = 1$.

\begin{prop} \label{inducción}
Let $\mathfrak{D}= \sum_{j=0}^{n}\partial^j  F_{j} \in \mathcal{D}(W)$ be a differential operator commuting with $D$. 
Then    
\begin{enumerate}
    \item [i) ] The leading coefficient of $\mathfrak D$ is $F_n=\begin{pmatrix}
    \alpha & a(\alpha-\beta)x \\ 0 &\beta
    \end{pmatrix}$, for some $\alpha, \beta \in \mathbb C$.
      
    \item [ii) ] For all $1 \leq k \leq n$,  we have
    $$ F_{n-k}= \begin{pmatrix}
    \tfrac{(-1)^{k}}{k!}\alpha\,[\![n]\!]_{k}\, x^{k} + h_{n-k} & r_{n-k} \\ 0 &\tfrac{(-1)^{k}}{k!}\beta\, [\![n]\!]_{k} \,x^{k} + g_{n-k}
    \end{pmatrix}, $$
    
 \noindent  where $h_{n-k},g_{n-k} \in \mathbb{C}[x]$ with $\operatorname{deg}(h_{n-k}) , \operatorname{deg}(g_{n-k}) \leq k-1$. 
    \end{enumerate}
\end{prop}
\begin{proof}

From Proposition \ref{(2,1)} we have that the coefficients of a differential operator in $\mathcal D(W)$ are upper triangular matrices, that is
$$F_{n-k} = \begin{pmatrix} p_{n-k} & r_{n-k} \\ 0 & q_{n-k} \end{pmatrix}.$$

The statement in i) has already been proven in \eqref{k=0}. 
 To see ii), we proceed by induction on $k\geq 1$. 
For $k = 1$, the statement is true from  \eqref{k=1}. 
Assume that the statement of the proposition is true for some $j$, i.e.
   \begin{equation}
       p_{n-j} =\tfrac{(-1)^{j}}{j!}\alpha\,[\![n]\!]_{j}\, x^{j} + h_{n-j},\qquad
      q_{n-j} = \tfrac{(-1)^{j}}{j!}\beta\, [\![n]\!]_{j} \,x^{j} + g_{n-j}. 
    \end{equation}
From \eqref{dD=Dd}, with $k=n-j$, we get
    \begin{equation}\label{paso inductivo}
        \begin{split}
         2F_{n-j-1}'+&F_{n-j-1}(Ax+B)-(Ax+B)F_{n-j-1} \\ &= (n-j)AF_{n-j}+CF_{n-j}-F_{n-j}C -F_{n-j}'(Ax+B) -F_{n-j}''           .
         \end{split}
    \end{equation}
    The left-hand side of \eqref{paso inductivo} is 
    \begin{equation}\label{pn-j-1}
        \begin{pmatrix}
        2p_{n-j-1}' & 2r_{n-j-1}'\\ 0&q_{n-j-1}'
        \end{pmatrix}
        + \begin{pmatrix}
        0&  -2br_{n-j-1}-(2abx-2a)(p_{n-j-1}-q_{n-j-1})\\0&0
        \end{pmatrix},
    \end{equation}
    and the right-hand side is 
    \begin{equation}\label{pn-j}
    \begin{split}
        -(n-j)&\begin{pmatrix}
        2p_{n-j}& 2ab q_{n-j}+2r_{n-j}\\0& 2q_{n-j}
        \end{pmatrix}+ \begin{pmatrix}
                0&-2r_{n-j}\\ 0&0
        \end{pmatrix}-\begin{pmatrix}
        p''_{n-j} &r''_{n-j}\\ 0&q''_{n-j}
        \end{pmatrix}\\
         & \quad +\begin{pmatrix}
        (2x-2b)p'_{n-j} & 2xr'_{n-j}+(2abx-2a)p_{n_j}'\\0 & 2xq_{n-j}'
        \end{pmatrix}.
    \end{split}
    \end{equation}
     By comparing the entries $(1,1)$ of the matrices  
     in \eqref{pn-j-1} and \eqref{pn-j},  and using the  inductive hypothesis,
    we obtain  
    \begin{align*}
        p_{n-j-1}'
        &= - \tfrac{(-1)^{j}}{j!} (n-2j) \, 
        [\! [n] \!]_{j}  \alpha\, x^{j} -(n-j)h_{n-j} + xh_{n-j}' - bp_{n-j}'- \tfrac 1{2}{p_{n-j}''}          .
       \end{align*}
 The right-hand side in the above equation is a polynomial of degree $\leq j$. Therefore 
    $$  p_{n-j-1}  = \frac{(-1)^{j+1}}{(j+1)!}[\! [n] \!]_{j+1} \alpha x^{j+1} + h_{n-j-1}, $$
where   $h_{n-j-1}=\int  \big (xh_{n-j}' -\frac1{2} {p_{n-j}''} -bp_{n-j}'- (n-j)h_{n-j}\big ) \, dx$. 
         By proceeding in the same way with the entry $(2,2)$,  we obtain 
    \begin{equation*}
            q_{n-j-1}  = \frac{(-1)^{j+1}}{(j+1)!}[\! [n] \!]_{j+1} \beta x^{j+1} + g_{n-j-1},
    \end{equation*}
    where $g_{n-j-1}=\int \big(- (n-j)g_{n-j} + xg_{n-j}' -\frac{q_{n-j}''}{2}\big) \, dx $.
\end{proof}

 \smallskip 
 
\begin{remark}
For $n=2m$ and $k>m$, the coefficient $[\![n]\!]_k=0$ and Proposition \ref{inducción}, does not give any new information because we already know that $\deg(F_{n-k})\leq n-k< k$.
\end{remark}
 
 \smallskip
 
\begin{prop} \label{caso impar}
    There are no-odd order operators $\mathfrak{D} \in \mathcal{D}(W)$ that commute with $D$.
\end{prop}
\begin{proof}
 Suppose that
  $\mathfrak{D} = \sum_{j=0}^{n}\partial^{j}F_{j} \in \mathcal{D}(W)$ is of odd order and it  
 commutes with $D$. Let us say  $n=2m-1$ and $F_n\neq 0$. 
 Recall that for any $\mathfrak D\in \mathcal{D} (W)$ we have $\deg(F_j)\leq j$, for all $0\leq j\leq n$. 
 
 From Proposition \ref{inducción} ii) with $k=m$, we get 
 $$F_{m-1} = \begin{pmatrix}
   \tfrac{(-1)^{m}}{m!}\alpha\,[\![2m-1]\!]_{m}\, x^{m} + h_{m-1} & r_{m-1} \\ 
        0&  \tfrac{(-1)^{m}}{m!}\beta\, [\![2m-1]\!]_{m} \,x^{m} + g_{m-1} \end{pmatrix},
        $$
     where $h_{m-1},g_{m-1}$ are polynomials  with $\deg(h_{m-1}) , \, \deg(g_{m-1}) \leq m-1$. 
   Since $[\![2m-1]\!]_{m}\neq0$, we have that $\alpha=\beta=0$. Therefore, from  Proposition \ref{inducción}, it follows that the leading coefficient of $\mathfrak D$ is 
   $F_n=0$, which is a contradiction.
\end{proof}

\begin{prop} \label{caso par}
    Let $\mathfrak{D}=\displaystyle \sum_{j=0}^{n}\partial^{j}F_{j} $ a differential operator in $ \mathcal D(W)$. 
    If $\mathfrak D $ commutes with $D$, then the leading coefficient of $\mathfrak D$ is scalar, i.e. $F_{n} = \alpha I$ for some $\alpha \in \mathbb{C}$.
\end{prop}
\begin{proof} 
Let $\mathfrak{D} \in \mathcal D(W)$ be a differential operator commuting with $D$. The differential operator $\mathfrak D^\dagger$ also commutes with $D$ because $D$ is symmetric, 
  and ${}^\dagger $ is an involution in $\mathcal D(W)$.  
  The operators 
   $\mathfrak D_1= \mathfrak{D}+\mathfrak{D}^{\dagger}$ and $\mathfrak D_2= i\mathfrak{D}-i\mathfrak{D}^{\dagger}$ are symmetric operators  commuting  with $D$ and $\mathfrak D=
  \frac 12\mathfrak D_1-i\frac 12 \mathfrak D_2$.
Therefore we can assume that  $\mathfrak D$ is a symmetric operator. 

    We write $F_{j} = \begin{pmatrix} p_{j} && r_{j} \\ 0 && q_{j} \end{pmatrix}$, for some polynomials $p_{j},r_{j},q_{j} \in \mathbb{C}[x]$ of degree less than or equal to 
    $j$, and $n = 2m$. From Proposition \ref{inducción}, with $k=m$, we have that 
    \begin{equation*}
        F_{m} = \begin{pmatrix} \frac{(-1)^{m}}{m!}[\![n]\!]_{m}\alpha x^{m} + h_{m} & r_{m} \\ 0 & \frac{(-1)^{m}}{m!}[\![n]\!]_{m}\beta x^{m} + g_{m} \end{pmatrix}
    \end{equation*}
   and from \eqref{r} 
    \begin{equation*}
    \begin{split}
        r_{m}  
          &= ax(q_{m}-p_{m}) + a \sum_{j=0}^{m-1}(-1)^{2m-j} \binom{2m -j }{ m} (m-j)e^{x^{2}}\,(e^{-x^{2}}\overline{q}_{2m-j})^{(m-j-1)}. 
    \end{split}
    \end{equation*}
    
    Given a polynomial $f$, the function $e^{x^{2}}\,(e^{-x^{2}} f(x) )^{(j)}$ is a polynomial of degree equal to $\deg(f)+j$.
    From Proposition \ref{inducción} ii), we have that $\operatorname{deg}(q_{2m-j}) = j$, for all $0\leq  j \leq m$. 
    Thus $\deg\big( (e^{-x^{2}}\overline{q}_{2m-j})^{(m-j-1)}e^{x^{2}}\big)=m-1$.

    If $\alpha\neq \beta$, then the polynomial $q_m-p_m$ is of degree $m$, and  $\deg(r_m)=m+1$, which is a contradiction. 
    Therefore,  $\alpha=\beta$, 
    $p_n=q_n=\alpha$,  and $ r_n=0$, which 
   concludes the proof of the proposition.
\end{proof}

\smallskip

Finally, we obtain that the centralizer $\mathcal Z_{\mathcal D(W)}(D)$ of $D$ in $\mathcal{D}(W)$, is a polynomial algebra in $D$.

\begin{thm} \label{R}
    Let $\mathfrak{D} \in \mathcal{D}(W)$ be a differential operator commuting with  D. Then $\mathfrak D$ is a polynomial in $D$. 
\end{thm}
\begin{proof} 
We proceed by induction on $m$, where $\operatorname{ord}(\mathfrak{D}) = 2m$.
Assuming that the proposition is true for differential operators of order $\leq 2(m-1)$, let $\mathfrak D=\sum_{j=0}^{2m}\partial^{j}F_{j}$ be an operator of order $2m$.
From Proposition \ref{caso par}, we have that $F_{2m} = \alpha I$ for some $\alpha \in \mathbb{C}$.
The differential operator $\mathfrak{B} = \mathfrak{D} - \alpha D^{m}$ commutes with $D$ and is of order less than or equal to $2(m-1)$ because there are no operators of odd order in the algebra $\mathcal{D}(W)$.
Thus, by the inductive hypothesis, we have that $\mathfrak{B}$ is a polynomial in $D$, and then $\mathfrak{D} \in \mathbb{C}[D]$. 
\end{proof}

The rest of this section will aim to prove that the entire algebra $\mathcal D(W) $ is a polynomial algebra in the differential operator $D$, in other words, that the centralizer of $D$ is the entire algebra $\mathcal D(W)$. 
Let $\mathfrak{D} = \sum_{j=0}^{s}\partial^{j}F_{j} \in \mathcal D(W)$, and let 
  $\{P_n(x)\}_{n\in \mathbb{N}_0}$ be the sequence of monic orthogonal polynomials for $W$. We have that
$$ P_n\mathfrak{D}=\Lambda_n P_n, \qquad \text{for all } n\in\mathbb{N}_0.$$
The eigenvalues $\Lambda_n= \Lambda_n(\mathfrak D)$ are determined in terms of the leading coefficients of the polynomials $F_i$. Explicitly, if 
 $F_i(x)=\sum_{j=0}^i x^j F_j^i$, 
 then the eigenvalues are obtained by
  $ \Lambda_n(\mathfrak D)=\sum_{i=0}^s [n]_i F_i^i,$ for all 
  $ n\geq 0$, where $[n]_i=n(n-1)\cdots (n-i+1)$.  
Hence, for any $\mathfrak{D}\in \mathcal D(W)$, the map $$n\longrightarrow \Lambda_n(\mathfrak D)$$ is a matrix-valued polynomial function of degree less or equal to $\operatorname{ord}(D)$. Moreover, from Proposition   \ref{(2,1)}, these eigenvalues are triangular matrices, let say $$\Lambda_{n}(\mathfrak{D}) = \begin{pmatrix} p(n) && r(n) \\ 0 && q(n) \end{pmatrix},$$ for some $p,q,r \in \mathbb{C}[x].$

\smallskip
\begin{prop}\label{autovalor}
A differential operator 
$\mathfrak{D}\in \mathcal{D}(W)$ commutes with $D$ if and only  if $$r(n)=a\,b\,n\,\big(p(n)-q(n)\big),$$
i.e. the eigenvalues of $\mathfrak D$ are of the form $$\Lambda_{n}(\mathfrak{D}) = \begin{pmatrix}p(n) && abn\,(p(n)-q(n)) \\ 0 && q(n) \end{pmatrix}.$$ 
\end{prop}
\begin{proof} 
The sequence of representations $\{ \Lambda_n\}_{n\in \mathbb{N}_0}$ separates points of the the algebra $\mathcal D(W)$. Hence,  $\mathfrak D$ commutes with $D$ if and only if $\Lambda_{n}(D)\Lambda_{n}(\mathfrak{D})=\Lambda_{n}(\mathfrak{D})\Lambda_{n}(D)$, for all $n\geq 0$.
The eigenvalues of $D$  and $\mathfrak D$ are given by
$$\Lambda_{n}(D) = \begin{pmatrix} -2n-2 && -2abn \\ 0 && -2n \end{pmatrix} \quad \text{ and } \quad 
\Lambda_{n}(\mathfrak{D}) = \begin{pmatrix} p(n) && r(n) \\ 0 && q(n) \end{pmatrix},$$
then we  have that 
\begin{equation*}
    \Lambda_{n}(D)\Lambda_{n}(\mathfrak{D})-\Lambda_{n}(\mathfrak{D})\Lambda_{n}(D) = \begin{pmatrix} 0 && -2r(n) + 2 abn\big (p(n)-q(n) \big) \\ 0 && 0 \end{pmatrix},
\end{equation*} 
which completes the proof of the proposition. 
\end{proof}

\smallskip

\begin{prop}\label{center}
The algebra $\mathcal D(W)$ has a non-trivial center. 
\end{prop}
\begin{proof}
From Theorem 4.13 in \cite{CY18}, the algebra $\mathcal{D}(W)$ is a finitely generated module over its center $\mathcal{Z}(W)$.
If $\mathcal{Z}(W) = \mathbb{C}I$ and $\mathfrak{D}_{1},\ldots,\mathfrak{D}_{\ell}\in \mathcal{D}(W)$ are such generators of $\mathcal D(W)$, then any operator in $\mathcal D(W)$ is a linear combination of them.
However, the operator $\sum_{k=1}^{\ell}w_{k}\mathfrak{D}_{k}$ is of order at most $M= \operatorname{max}\{\operatorname{ord}(\mathfrak{D}_{1}),\ldots,\operatorname{ord}(\mathfrak{D}_{\ell}) \}$, which leads to a contradiction.
\end{proof}

\begin{thm}\label{polynomial algebra}
The algebra $\mathcal{D}(W)$ coincides with the centralizer of $D$ in $\mathcal{D}(W)$. 
In particular, it is a polynomial algebra in the differential operator $D$. 
\end{thm}
\begin{proof}

Let $\mathfrak{D} \in \mathcal{D}(W)$. The sequence of eigenvalues of the monic orthogonal polynomials is given by $\Lambda_{n}(\mathfrak{D}) = \begin{pmatrix} p(n) && r(n) \\ 0 && q(n) \end{pmatrix}$, for some $p,q,r \in \mathbb{C}[x]$.
From Propositions \ref{center} and \ref{autovalor}, there exists a differential operator $E \in \mathcal{Z}(W)$ with $\operatorname{ord}(E)>0$, 
and $$\Lambda_{n}(E) =\begin{pmatrix}s(n) && ab\,n\,(s(n)-t(n)) \\0 && t(n) \end{pmatrix},$$
for some $s,t $ polynomials of degree at most $\operatorname{Ord}(E)$.
In particular, $s(n)-t(n)\neq 0$ almost everywhere.
Since $\Lambda_n(E) \Lambda_{n}(\mathfrak{D})=\Lambda_{n}(\mathfrak{D})\Lambda_{n}(E)$,  we have that 
$$\Big ( r(n) - abn\,(p(n)-q(n)) \Big )(s(n)-t(n))=0.$$ 
Therefore, $r(n) = abn\,(p(n)-q(n)$, for all  $n$, and  
the theorem follows from Proposition \ref{autovalor}.  
\end{proof}

\smallskip 
Recall that the algebra $\mathcal{D}(W)$ is  full if 
there exist nonzero $W$-symmetric operators $\mathfrak{D}_{1},\ldots,\mathfrak{D}_{N}$ in $\mathcal{D}(W)$,  such that  $$\mathfrak{D}_{i}\mathfrak{D}_{j} = 0 \text{ for $i\not= j\quad $ with } \quad \mathfrak{D}_{1} + \cdots + \mathfrak{D}_{N} \in \mathcal{Z}(W)$$ 
which is not a zero divisor. 
As a consequence of Theorem \ref{polynomial algebra}, we obtain the following result.

\begin{thm}\label{not full} The algebra $\mathcal D(W)$ is not a full algebra.
\end{thm}
\begin{proof}
Let $\mathfrak{D}_{1},\mathfrak{D}_{2} $ be nonzero $W$-symmetric operators in $\mathcal{D}(W)$.  
From Theorem \ref{polynomial algebra}, we have that $\mathfrak{D}_{1} = \sum_{j =0}^{n}\alpha_{j}D^{j}$, $\mathfrak{D}_{2} = \sum_{k=0}^{m} \beta_{k}D^{k}$ are polynomials in the differential operator $D$.
Thus we have that the leading coefficient of $\mathfrak{D}_{1}\mathfrak{D}_{2}$ is $\alpha_{n}\beta_{m}I \not= 0$, and therefore $\mathfrak{D}_{1}\mathfrak{D}_{2} \not=0$. 
\end{proof}

As a direct consequence of Theorem \ref{clas1} and Theorem \ref{not full}, we conclude that  

\begin{thm}\label{not Darboux2}
The weight matrix 
$$W(x)=e^{-x^{2}}\begin{pmatrix} e^{2bx} + a^{2}x^{2} && ax \\ ax && 1 \end{pmatrix}\qquad  (a,b\neq 0)$$
is not a noncommutative bispectral Darboux transformation of any direct sum of classical weights.
\end{thm}

\section{The right Fourier Algebra of $W(x)$ }\label{RFA}

 In this section, we study the right Fourier algebra of the weight $W(x)$. 
Even when  $W$ is not a Darboux transformation of a direct sum of classical weights, it is closely related to it, as we explain below.  
We will take full advantage of this fact to obtain an explicit description of the right Fourier algebra associated with $W$. 

By introducing  the diagonal weight 
\begin{equation}\label{Wtilde}
    \widetilde{W}(x)  = \begin{pmatrix} e^{-x^2+2bx} && 0 \\ 0 && e^{-x^2} \end{pmatrix},
   \end{equation}
 we get the following factorization 
\begin{equation}\label{WWtilde}
W(x) = T(x)\widetilde{W}(x)T(x)^{\ast} \qquad \text { with } \quad 
T(x) = \begin{pmatrix}1 && ax \\ 0 && 1 \end{pmatrix}.
\end{equation}

From  \eqref{fourier algebra} and  Remark \ref{adj}, a  differential operator $\mathfrak{D} \in \operatorname{Mat}_{2}(\Omega[x])$ belongs to
$ \mathcal{F}_{R}({W})$ if and only if its $W$-formal adjoint 
$\mathfrak{D}^{\dagger}  = {W}(x) \mathfrak{D}^{\ast}{W}(x)^{-1} $ has  polynomial coefficients. 
Similarly, $\mathfrak{\widetilde D} \in \operatorname{Mat}_{2}(\Omega[x])$ belongs to
$ \mathcal{F}_{R}({\widetilde W})$ if and only if 
$\mathfrak{\widetilde D}^{\dagger}  = {\widetilde W}(x) \mathfrak{\widetilde D}^{\ast}{\widetilde W}(x)^{-1} \in \operatorname{Mat}_{2}(\Omega[x])$.

\begin{prop}\label{relationFourier} The right Fourier algebras of $W$ and $\widetilde W$ satisfy the relation  
\begin{equation*}
    \mathcal{F}_{R}(W) = T(x)\mathcal{F}_{R}(\widetilde{W})T^{-1}(x).
\end{equation*}
\end{prop}
\begin{proof} 
    For any  $\mathfrak{D}\in \mathcal{F}_{R}(W)$      
the differential operators      
$\mathfrak{D}^{\dagger}   =  W(x)\mathfrak{D}^{\ast}W^{-1}(x) $ and    
$E = T^{-1}(x)\mathfrak{D}T(x)$ have also polynomial coefficients because $T(x)$ and $T^{-1}(x)$ are  polynomials. We also get 
$$T^{-1}(x)\, \mathfrak{D}^\dagger\, T(x)=\widetilde{W}(x) E^{\ast}\,\widetilde{W}^{-1}(x)=E^\dagger. $$
Therefore, $ E\in \mathcal{F}_{R}(\widetilde W)$ and $\mathfrak D\in T(x)\mathcal{F}_{R}(\widetilde{W})T^{-1}(x)$.
\end{proof}

\begin{remark}
We get that $$\mathcal D(W)\subsetneq T(x)\mathcal D(\widetilde W)T^{-1}(x).$$
In fact, the generator $D$ of the algebra $\mathcal D(W)$ can be factorized as $D= T(x)\widetilde D T^{-1}(x)$ with 
$$\widetilde{D} = \begin{pmatrix} \partial^{2} + \partial(-2(x-b)) - 2 && 0 \\ 0 && \partial^{2} + \partial(-2x) \end{pmatrix}\in \mathcal{D}(\widetilde{W}). $$
Observe that $\partial^{2} + \partial(-2x)$ is the classical Hermite operator. 
Then we have that $\widetilde{D}$
is a $\widetilde{W}$-symmetric differential operator 
in $\mathcal{D}(\widetilde{W})$. 
On the other hand the operator 
$$\widetilde{D_2} = \begin{pmatrix} \partial^{2} + \partial(-2(x-b))  && 0 \\ 0 && \partial^{2} + \partial(-2x) \end{pmatrix}\in \mathcal{D}(\widetilde{W}),$$
but 
$$T(x)\widetilde{D}_2 T(x)^{-1}= \partial^2 I+ \partial \, \begin{pmatrix} -2x+2b & -2abx+2a\\ 0& -2x
\end{pmatrix} +\begin{pmatrix}
    0&-2ax\\ 0&0
\end{pmatrix}\not \in \mathcal D (W).$$ 
   
\end{remark}

\smallskip

The following results give an explicit description of the right Fourier algebra of the weight $W$.

\begin{thm}\label{Fourier-W}
The right Fourier algebra of ${W}$ is
$$\mathcal{F}_{R}(W) =\left \{ \sum_{j=0}^{n}  \partial^{j} \begin{pmatrix} p_{j}(x) && ax(q_{j}(x)-p_{j}(x)) \\ 0 && q_{j}(x) \end{pmatrix} + \partial^{j-1}\begin{pmatrix} 0 && jaq_{j}(x) \\ 0 && 0 \end{pmatrix} \, : \,  p_{j},q_{j} \in \mathbb{C}[x], \, n\geq 0 \right \}.$$
\end{thm}

\begin{proof} 
From Proposition \ref{relationFourier} and \eqref{WWtilde} it is enough to prove that 
the right Fourier algebra of $\widetilde{W}$ is $$\mathcal{F}_{R}(\widetilde{W}) =\left \{ \sum_{j=0}^{n} \partial^{j} F_{j}(x) \, :
\, 
F_{j}(x) = \begin{pmatrix} p_{j}(x) && 0 \\ 0 && q_{j}(x) \end{pmatrix}, p_{j},q_{j} \in \mathbb{C}[x], \, n\geq 0 \right \}.$$
In fact, a differential operator $\mathfrak{D}=\sum_{j=0}^{n} \partial^{j}F_{j}(x)
$ with polynomial coefficients belongs to
$ \mathcal{F}_{R}(\widetilde{W})$ if and only if 
$\mathfrak{D}^{\dagger}  = \widetilde{W}(x) \mathfrak{D}^{\ast}\widetilde{W}(x)^{-1} 
$ has also polynomial coefficients. We have that 
\begin{align*}
    \mathfrak{D}^{\dagger}& 
=\sum_{j=0}^{n} (-1)^{j} \sum_{k=0}^{j}\binom{j}{k}\, \partial ^{k}\big (\widetilde{W}(x)F_{j}(x)^{\ast}\big )^{(k)}\widetilde{W}(x)^{-1}
= \sum_{k=0}^{n}\partial^{k}G_{k}(x), 
\end{align*} where 
\begin{equation}\label{GG}
    G_k=  \sum_{j = 0}^{n-k}(-1)^{n-j} \binom{n-j}{k}
    (\widetilde{W}(x)F_{n-j}(x)^*)^{(n-k-j)} \widetilde{W}(x)^{-1} , \qquad \text{for } 0\leq k\leq n.  
 \end{equation}
If $F_j=\begin{pmatrix} p_{j}(x) && 0 \\ 0 && q_{j}(x) \end{pmatrix} $ is a polynomial, then 
\begin{equation}\label{cof}
(\widetilde{W}(x)F_{j}(x)^{\ast})^{(k)}\widetilde{W}(x)^{-1} = \begin{pmatrix} (e^{-x^{2}+2bx}\overline{p}_{j}(x))^{(k)}e^{x^{2}-2bx} && 0 \\ 0 && (e^{-x^{2}}\overline{q}_{j}(x))^{(k)}e^{x^{2}}\end{pmatrix}
\end{equation} 
is also a polynomial. Thus, $\mathfrak{D}^{\dagger} \in \operatorname{Mat}_{2}(\Omega[x])$ which implies that $\mathfrak{D}\in \mathcal{F}_{R}(\widetilde{W})$.

On the other hand, for  $\mathfrak{D}= \sum_{j=0}^{n} \partial^{j}F_{j}(x)\in \mathcal{F}_{R}(\widetilde{W})$ we will see that all  $F_j$ are diagonal matrices. We have that  
$\mathfrak{D}^{\dagger} = \sum_{j=0}^{n} \partial^{j}G_{j}(x)$
has polynomial coefficients. 
If $F_{j}(x) = \begin{pmatrix} p_{j}(x) && r_{j}(x) \\ g_{j}(x) && q_{j}(x) \end{pmatrix} \in \operatorname{Mat}_{2}(\mathbb{C}[x])$,
by $\eqref{GG}$ we get 
$$G_{n}(x) =(-1)^{n} \widetilde{W}(x)F_{n}(x)^{\ast}\widetilde{W}(x)^{-1} =(-1)^{n} \begin{pmatrix} \overline{p}_{n}(x) && e^{2bx}\overline{g}_{n}(x) \\ e^{-2bx}\overline{r}_{n}(x) && \overline{q}_{n}(x)  \end{pmatrix}.$$
Hence 
$g_{n}(x)=0$ and $r_{n}(x)=0$.  
By induction on $k$ we can assume that $g_{n-j}(x)=0$ and $r_{n-j}(x)=0$, i.e. $F_{n-j}$ are diagonal, for all $0 \leq j \leq k$. 
By \eqref{GG} we have that 
$$G_{n-k-1}  = \sum_{j=0}^{k+1}(-1)^{n-j}\binom{n-j}{n-k-1}(\widetilde{W}(x)F_{n-j}(x)^{\ast})^{(k+1-j)}\widetilde{W}(x)^{-1}.$$
For $0 \leq j \leq k$,   we get 
$$\sum_{j=0}^{k}(-1)^{n-j}\binom{n-j}{n-k-1}(\widetilde{W}(x)F_{n-j}(x)^{\ast})^{(k+1-j)}\widetilde{W}(x)^{-1} = \begin{pmatrix} \alpha(x) && 0 \\ 0 && \beta(x) \end{pmatrix},$$
for some polynomials  $\alpha,\beta$. 
Thus 
\begin{equation*}
    \begin{split}
        G_{n-k-1} &  = \begin{pmatrix} \alpha(x) + (-1)^{n-k-1}\overline{p}_{n-k-1}(x) && (-1)^{n-k-1}e^{2bx}\overline{g}_{n-k-1}(x) \\ (-1)^{n-k-1}e^{-2bx}\overline{r}_{n-k-1}(x) && \beta(x) + (-1)^{n-k-1}\overline{q}_{n-k-1}(x) \end{pmatrix}.
    \end{split}
\end{equation*}
is a matrix polynomial which implies that 
$g_{n-k-1}(x)=0$ and $r_{n-k-1}(x)=0$, and this concludes the proof.
\end{proof}

\begin{prop}
  Let $\widetilde W$ be the diagonal weight matrix given in \eqref{Wtilde}.  
  The algebra $\mathcal D(\widetilde W)$  is a commutative full algebra.
  \end{prop}

\begin{proof} The algebra $\mathcal D(\widetilde W) $ of any diagonal weight matrix $\widetilde{ W}$ is always a full algebra. In fact,
    $\mathcal D_1=\text{diag} (1,0)$ and   $\mathcal D_2=\text{diag} (0,1)$ are $\widetilde W$-symmetric operators in $\mathcal D
    (\widetilde W)$ satisfying $\mathcal D_1\mathcal D_2=0=\mathcal D_2\mathcal D_1$ and $\mathcal D_1+\mathcal D_2$ is a central element. 
    
    On the other hand, $\mathcal D(\widetilde W)\subset \mathcal F_R(\widetilde W)$, which implies that  $\mathcal D(\widetilde W) $ is a diagonal algebra. Thus $\mathcal D(\widetilde W)= \mathcal D(e^{-x^2+2bx})\oplus \mathcal D(e^{-x^2})$.   
    Therefore, $\mathcal D(\widetilde{W})$ is a commutative algebra. 
\end{proof}

\begin{remark}
   Observe that $\widetilde W$ is a noncommutative bispectral Darboux transformation of scalar weights, and the algebra $\mathcal D(\widetilde W)$ is commutative (cf. Theorem 1.5 in \cite{CY18}).
\end{remark} 

\section{Orthogonal polynomials with respect to $W$.}
\label{mop}

In this section, we give the explicit expressions of a 
sequence of orthogonal polynomials for the weight 
$$W(x)=e^{-x^{2}}\begin{pmatrix} e^{2bx} + a^{2}x^{2} && ax \\ ax && 1 \end{pmatrix},$$ 
and we also give the three-term recursion relation satisfied by them.

\smallskip
We denote $H_n(x)$ the $n$-th monic Hermite polynomial, given by 
$H_{n}(x) = \frac{(-1)^{n}}{2^{n}}(e^{-x^{2}})^{(n)}e^{x^{2}}.$ We observe that $H_{n}(x-b)$ is the sequence of the monic orthogonal polynomials for the weight $\tilde w(x)=e^{-x^{2}+2bx}$.

\begin{prop}\label{orthpoly} The polynomials
$$Q(x,n) = \begin{pmatrix}H_{n}(x-b) && aH_{n+1}(x)-axH_{n}(x-b) \\ -anH_{n-1}(x-b) && 2e^{b^{2}}H_{n}(x)+a^{2}nxH_{n-1}(x-b)\end{pmatrix}$$ are orthogonal with respect to $W(x)$. Moreover,  the leading coefficient of $Q(x,n)$ is given by 
\begin{equation}\label{Mn}
    M_{n} = \begin{pmatrix} 1 && nab \\ 0 && 2e^{b^{2}}+a^{2}n \end{pmatrix}.
\end{equation}
\end{prop}

\begin{proof}
The weight $W(x)$ can be  factorized as  
$W(x) = T(x)\widetilde{W}(x)T(x)^{\ast}$, where 
 $$ T(x)=\begin{pmatrix} 1 && ax \\ 0 && 1 \end{pmatrix} \quad \text{ and } \quad  \widetilde{W}(x)=\begin{pmatrix} \tilde w(x) && 0 \\ 0 && w(x) \end{pmatrix},$$
 and $w(x)=e^{-x^2}$, 
 $\tilde w(x)=e^{-b^2}w(x-b)$.  Thus
$$Q(x,n)T(x) 
= \begin{pmatrix} H_{n}(x-b) && aH_{n+1}(x) \\ -anH_{n-1}(x-b) && 2e^{b^{2}}H_{n}(x) \end{pmatrix}.$$
We have that the sequence $\{Q(x,n)\}_{n}$ is orthogonal with respect to $W(x)$ if and only if $\{Q(x,n)T\}_n$ is orthogonal with respect to $\widetilde{W}$.

For $n\neq m$, we compute 
\begin{align*}
 \langle & Q(x,n)T(x),  Q(m,x)T(x)\rangle_{\widetilde{W}}=\\
 & \quad  = \int_{-\infty}^{\infty}\begin{pmatrix} H_{n}(x-b) && aH_{n+1}(x) \\ -anH_{n-1}(x-b) && 2e^{b^{2}}H_{n}(x) \end{pmatrix}\widetilde{W}(x)\begin{pmatrix} H_{m}(x-b) && -amH_{m-1}(x-b) \\  aH_{m+1}(x) && 2e^{b^{2}}H_{m}(x) \end{pmatrix}dx 
 \\ & \quad = \begin{pmatrix} k_{11} && k_{12} \\ k_{21} && k_{22} \end{pmatrix}.
\end{align*}

 \ 

\noindent 
To see that 
$$k_{11}=\langle H_{n}(x-b),H_{m}(x-b) \rangle_{\widetilde{w}} + a^{2} \langle H_{n+1}(x),H_{m+1}(x) \rangle_{w}=0, $$ and 
$$ k_{22} = nm a^2\langle H_{n-1}(x-b), H_{m-1}(x-b) \rangle_{\widetilde{w}} + 4e^{2b^{2}} \langle H_{n}(x), H_{m}(x) \rangle_{w}= 0, $$
we use that $\{H_n(x)\}_n$ and $\{H_n(x-b)\}_n$ are orthogonal with respect to $w$ and $\tilde w$ respectively.

Now we compute  $$k_{12}=-am \langle H_{n}(x-b), H_{m-1}(x-b) \rangle_{\widetilde{w}} + 2ae^{b} \langle H_{n+1}(x),H_{m}(x) \rangle_{w}.$$ 
If $m \not= n+1$, then  $k_{12}=0$.
If $m = n + 1$, we use that 
$\|H_{n+1}\|^{2}_{w} =\frac{(n+1)}{2}\|H_{n}\|^{2}_{w}$, and 
$$\langle H_{n}(x-b), H_{n}(x-b)\rangle_{\widetilde{w}}  = \int_{-\infty}^{\infty}H_{n}(x-b)e^{-x^{2}+2bx}H_{n}(x-b)dx 
= e^{b^{2}}\|H_{n}\|^{2}_{w}
  $$
to obtain that $k_{12}=0$.
We proceed similarly to prove that 
$$ k_{21} = -an \langle H_{n-1}(x-b), H_{m}(x-b) \rangle_{\widetilde{w}} + 2ae^{b} \langle H_{n}(x),H_{m+1}(x) \rangle_{w}=0,$$
which concludes the proof that $Q(x,n)$ is an orthogonal sequence of polynomials for the weight $W$. The last assertion in the statement follows easily.
\end{proof}

\begin{prop}\label{ttrrQ}
The matrix orthogonal polynomials $Q(x,n)$ defined in Proposition \ref{orthpoly} satisfy the three-term recursion relation 
$$Q(x,n)x = \widetilde{A}_{n}Q(n+1,x) + \widetilde{B}_{n}Q(x,n)+\widetilde{C}_{n}Q(n-1,x)$$
where 
$$\widetilde{A}_{n}= \begin{pmatrix}1 && - \frac{ab}{2e^{b^{2}}+(n+1)a^{2}} \\ 0 && \frac{2e^{b^{2}}+na^{2}}{2e^{b^{2}}+(n+1)a^{2}} \end{pmatrix},  
\widetilde{B}_{n} = \begin{pmatrix} \frac{2e^{b^{2}}b}{2e^{b^{2}}+(n+1)a^{2}} &&  \frac{a}{2(2e^{b^{2}}+na^{2})} \\ \frac{2e^{b^{2}}a}{2e^{b^{2}}+(n+1)a^{2}} && \frac{na^{2}b}{2e^{b^{2}}+na^{2}}\end{pmatrix}, 
\widetilde{C}_{n} = \begin{pmatrix} \frac{n}{2} \frac{(n+1)a^{2}+2e^{b^{2}}}{2e^{b^{2}}+na^{2}} && 0 \\ - \frac{2nabe^{b^{2}}}{2e^{b^{2}}+na^{2}} && \frac{n}{2} \end{pmatrix}.$$
\end{prop}
\begin{proof}
Let  $T = T(x) = \begin{pmatrix} 1 && ax \\ 0 && 1 \end{pmatrix}.$ Recall  that $Q(x,n)T = \begin{pmatrix} H_{n}(x-b) && aH_{n+1}(x) \\ -anH_{n-1}(x-b) && 2e^{b^{2}}H_{n}(x) \end{pmatrix}$. 

\noindent We compute

$$\widetilde{A}_{n}Q(n+1,x)T = \begin{pmatrix} H_{n+1}(x-b) + \frac{a^{2}b(n+1)}{a^{2}(n+1)+2e^{b^{2}}}H_{n}(x-b) && aH_{n+2}(x)- \frac{2abe^{b^{2}}}{a^{2}(n+1)+2e^{b^{2}}}H_{n+1}(x) \\ -\frac{(2e^{b^{2}}+na^{2})a(n+1)}{(n+1)a^{2}+2e^{b^{2}}}H_{n}(x-b)&& \frac{2(2e^{b^{2}}+na^{2})e^{b^{2}}}{(n+1)a^{2}+2e^{b^{2}}} H_{n+1}(x) \end{pmatrix},$$

$$\widetilde{B}_{n}Q(n.x)T = \begin{pmatrix}
\frac{2e^{b^{2}}b\,H_{n}(x-b)}{(n+1)a^{2}+2e^{b^{2}}} -\frac{a^{2}nH_{n-1}(x-b)}{2(2e^{b^{2}}+na^{2})} && \frac{2abe^{b^{2}}H_{n+1}(x)}{(n+1)a^{2}+2e^{b^{2}}} + \frac{ae^{b^{2}}H_{n}(x)}{2e^{b^{2}}+na^{2}} \\ \frac{2ae^{b^{2}}H_{n}(x-b)}{(n+1)a^{2}+2e^{b^{2}}}-\frac{n^{2}ba^{3}H_{n-1}(x-b)}{2e^{b^{2}}+na^{2}} && \frac{2a^{2}e^{b^{2}}H_{n+1}(x)}{(n+1)a^{2}+2e^{b^{2}}} + \frac{2na^{2}be^{b^{2}}H_{n}(x)}{2e^{b^{2}}+na^{2}} \end{pmatrix},$$

$$\widetilde{C}_{n}Q(n-1,x)T = \left(\begin{smallmatrix}\frac{n((n+1)a^{2}+2e^{b^{2}})
}{2(2e^{b^{2}}+na^{2})} H_{n-1}(x-b)&& \frac{n((n+1)a^{2}+2e^{b^{2}})a}{2(2e^{b^{2}}
+na^{2})} H_{n}(x)\\ -\frac{2anbe^{b^{2}}}{2e^{b^{2}}+na^{2}}H_{n-1}(x-b)
-\frac{n(n-1)a}{2}H_{n-2}(x-b) && -\frac{2na^{2}be^{b^{2}}}{2e^{b^{2}}+na^{2}}H_{n}(x) + ne^{b^{2}}H_{n-1}(x) \end{smallmatrix}\right).$$

Let $M=\widetilde{A}_{n}Q(n+1,x)T + \widetilde{B_{n}}Q(x,n)T + \widetilde{C}_{n}Q(n-1,x)T = \begin{pmatrix}m_{11} && m_{12} \\ m_{21} && m_{22} \end{pmatrix}$. We get
\begin{equation*}
    \begin{split}
         m_{11} & = H_{n+1}(x-b)+ bH_{n}(x-b)+\frac{n}{2}H_{n-1}(x-b), \\
         m_{12} & = a  \left ( \frac{(n+1)}2H_{n}(x)+H_{n+2}(x) \right ), \\  
         m_{21} &= -an \left ( H_{n}(x-b)+bH_{n-1}(x-b)+\frac{(n-1)}2H_{n-2}(x-b) \right ) , \\
         m_{22} &= e^{b^{2}}(2H_{n+1}(x) + nH_{n-1}(x)).
    \end{split}
\end{equation*}
By using the three-term recursion relation satisfied by the monic orthogonal Hermite polynomials 
\begin{equation*}
    H_{n}(x)x = H_{n+1}(x) + \frac{n}{2}H_{n-1}(x),
\end{equation*} 
we get that $M=x\,Q(x,n)T$. Thus, multiplying by $T^{-1}$ we complete the proof. 
\end{proof}

\smallskip
The monic orthogonal polynomials with respect to $W(x)$ are $P(x,n)= M_n^{-1}Q(x,n)$, with $M_{n}$ the matrix given in \eqref{Mn}. 
By taking 
$$A_{n} = M_{n}^{-1}\widetilde{A}_{n}M_{n+1}, \quad   
B_{n} = M_{n}^{-1}\widetilde{B}_{n}M_{n},\quad  C_{n} = M_{n}^{-1}\widetilde{C}_{n}M_{n-1}$$ we get 
$$P(x,n)x = A_n P(x,n+1) + B_{n} P(x,n)+ {C}_{n}P(x,n-1).$$
 Explicitly, we have $A_n=I$,
 \begin{align*}
     B_{n} &= \begin{pmatrix} \frac{4be^{2b^{2}}}{(na^{2}+a^{2}+2e^{b^{2}})(2e^{b^{2}}+na^{2})}&& \frac{a^{5}n(n+1)(2b^2n-1)+a^{3}2e^{b^{2}}(2b^{2}n^{2}-2n-1)-ae^{2b^{2}}(2b^{2}n+1)}{-2(na^{2}+a^{2}+2e^{b^{2}})(2e^{b^{2}}+na^{2})} \\ \frac{2e^{b^{2}}a}{(2e^{b^{2}}+na^{2})(na^{2}+a^{2}+2e^{b^{2}})} && \frac{na^{2}b(na^{2}+a^{2}+4e^{b^{2}})}{(na^{2}+a^{2}+2e^{b^{2}})(2e^{b^{2}}+na^{2})}\end{pmatrix},\\
     C_{n} &= \begin{pmatrix}\frac{n(a^{4}n(n+1)+a^{2}2e^{b^{2}}(2b^{2}n+2n+1)+4e^{2b^{2}})}{2(2e^{b^{2}}+na^{2})^{2}} && \frac{nab(a^{4}n(n-1)+a^{2}2e^{b^{2}}(2b^{2}n^{2}-2b^{2}n-1)-4e^{2b^{2}})}{2(2e^{b^{2}}+na^{2})^{2}} \\ \frac{-2nabe^{b^{2}}}{(2e^{b^{2}}+na^{2})^{2}} && \frac{-n(-a^{4}n(n-1)+a^{2}2e^{b^{2}}(2b^{2}n-2b^{2}-2n+1)-4e^{2b^{2}})}{2(2e^{b^{2}}+na^{2})^{2}} \end{pmatrix}.
 \end{align*}

\section{Higher dimensional examples}\label{HDEx}

Examples of weight matrices $W(x)$ which are solutions of the Matrix Bochner Problem and which are not obtained as noncommutative bispectral Darboux transformations of diagonal (classical) weights are also present in higher dimensions. 
With the help of computer software, we have found new families of such examples.  
For example, for   $N=3$,  and real parameters $a_1, a_2, b_1, b_2 \in \mathbb{R}-\{0\}$ with $b_{1} \not= b_{2}$, we have that 
 the second-order differential operator 
\begin{equation*}
    \begin{split}
       D & = \partial^{2} I + \partial \begin{pmatrix}  -2x &&  2a_{1}b_{1}x+2a_{1} && 0 \\  0 &&   2b_{1}-2x &&  0 \\  0 && (2a_{2}b_{1}-2a_{2}b_{2})x+2a_{2} && 2b_{2}-2x \end{pmatrix} + \begin{pmatrix}  0 &&  2a_{1}b_{1} && 0 \\ 0 && 2 && 0 \\ 0 && 2a_{2}b_{1} && 0 \end{pmatrix}.
    \end{split}
\end{equation*}
is a symmetric operator with respect to 
\begin{equation*}
    \begin{split}
        W(x) &= e^{-x^{2}}\begin{pmatrix}  
   1+a_{1}^{2}x^{2} e^{2b_{1}x} && a_{1}x e^{2b_{1}x}&& a_{1}a_{2}x^{2} e^{2b_{1}x}\\  a_{1}x e^{2b_{1}x} &&  e^{2b_{1}x} && a_{2}x e^{2b_{1}x} \\ a_{1}a_{2}x^{2} e^{2b_{1}x} && a_{2}x e^{2b_{1}x}&& a_{2}^{2}x^{2} e^{2b_{1}x}+e^{2b_{2}x} \end{pmatrix}. 
   \end{split}
\end{equation*}

By following the same reasoning as in the case of $N=2$, but with more intricate and lengthy computations, we can demonstrate that the algebra $\mathcal D(W)$ is a polynomial algebra in $D$. Therefore, $W$ is not a noncommutative bispectral Darboux transformation of a direct sum of classical weights.

We expect that these specific examples will contribute to a more comprehensive understanding of the solution to the Matrix Bochner Problem.

\section{Acknowledgments}
We express our gratitude to I. Zurrián, F. A. Gr\"unbaum and R. Miatello for their insightful conversations, encouragement, and assistance in reviewing the early versions of this paper.

This paper was partially supported by SeCyT-UNC, CONICET, PIP 1220150100356.

\end{document}